\documentclass[12pt]{article}

\usepackage{a4wide}
\usepackage{amssymb}
\usepackage{amsfonts}
\usepackage{amsmath}
\input xy
\xyoption{arrow} \xyoption{matrix}

\date{}

\newtheorem{proposition}{Proposition}[section]
\newtheorem{theorem}[proposition]{Theorem}
\newtheorem{lemma}[proposition]{Lemma}

\newtheorem{corollary}[proposition]{Corollary}

\def\der{\partial }

\def\nFM0{{\nu }_{F,M_0}}
\def\nFN0{{\nu }_{F,N_0}}
\def\nGN0{{\nu }_{G,N_0}}

\def\N0{ {\bf N}_0 }

\def\t{\otimes}
\def\g{\gamma}

\def\ra{\rightarrow}

\def\Xpm{X^{\pm }}

\def\s{\sigma}
\def\Z{\mathbb{Z}}

\def\l1{{\lambda}_1}

\def\a{\alpha}
\def\a0{ {\alpha }_0}
\def\a1{ {\alpha }_1}

\def\l{\lambda}

\def\nFGM0{{\nu }_{F,G,M_0}}

\def\nFN0{{\nu}_{F,N_0}}

\def\sm{{\sigma}^m}

\def\sm1{{\sigma}^{-1}}

\def\smtp1{{\sigma}^{-t+1}}

\def\S1{S^{-1}}

\def\Xpm1{X^{\pm 1}_1}

\def\sPM1{{\sigma }^{\pm 1}}
\def\sMP1{{\sigma }^{\mp 1 }}

\def\d{\delta}

\def\di{{\rm d.ind}}

\def\L{\Lambda}

\def\G{\Gamma}

\def\CD{{\cal D}}

\def\Ytm1{Y^{t-1}}
\def\Yim1{Y^{i-1}}

\def\Aut{{\rm Aut}}

\def\ker{ {\rm ker } }

\def\CJ{ {\cal J}}

\def\D{ \Delta }

\def\SL2Z{ {\rm SL}_2({\bf Z}) }

\def\th{ \theta }

\def\Gp1{ G^{1 , 1 } }
\def\P11{ P^{-1 , 1 } }
\def\Pp1{ P^{1 , 1 } }

\def\th{\theta}

\def\nCLsr{{}^\nu\kern-2pt {\cal L}^{\sigma , \rho  }}
\def\nP{{}^\nu \kern-2pt P}
\def\nL{{}^\nu\kern-2pt L}
\def\nLL{{}^\nu\kern-2pt \Lambda}
\def\nPsr{{}^\nu\kern-2pt P^{\sigma , \rho  }}
\def\nLsr{{}^\nu\kern-2pt L^{\sigma , \rho  }}
\def\nuCL{{}^\nu\kern-2pt  {\cal L}}
\def\nCLsr{{}^\nu\kern-2pt {\cal L}^{\sigma , \rho  }}
\def\nCL1m{{}^\nu\kern-2pt {\cal L}^{-1 , 1  }}
\def\x1nu{x^\frac{1}{\nu}}
\def\xm1nu{x^{-\frac{1}{\nu}}}

\def\ra{\rightarrow }

\def\nAM0{{\nu }_{{\cal A},M_0}}
\def\nAN0{{\nu }_{{\cal A},N_0}}

\def\CJ{ {\cal J }}

\def\det{ {\rm det }}

\def\bx{\overline{x}}

\def\gp{\mathfrak{p}}

\def\GL{{\rm GL}}
\def\SL{{\rm SL}}

\def\Spec{{\rm Spec}}

\def\di!{\frac{\der^i}{i!}}
\def\dik!{\frac{\der^k_i}{k!}}

\def\Fp{\mathbb{F}_p}

\def\N{\mathbb{N}}

\def\0{\overline{0}}
\def\1{\overline{1}}

\def\Ln1{\L_{n,\overline{1}}}

\def\a1{a_{\overline{1}}}

\def\S{\Sigma}

\def\vn1{\overrightarrow{n-1}}

\def\Aff{{\rm Aff}}

\def\im{{\rm im}}

\def\res{{\rm res}}

\begin{document}

\author{V. V. \  Bavula 
}

\title{The group of automorphisms of the first Weyl algebra in prime characteristic and the restriction map }

\maketitle
\begin{abstract}
Let $K$ be a {\em perfect} field of characteristic $p>0$,
$A_1:=K\langle x, \der \, | \, \der x- x\der =1\rangle$ be the
first Weyl algebra and  $Z:=K[X:=x^p, Y:=\der^p]$ be its centre.
It is proved that $(i)$ the restriction map $\res :\Aut_K(A_1)\ra
\Aut_K(Z), \; \s \mapsto \s|_Z$, is a monomorphism with $\im (\res
) = \G :=\{ \tau \in \Aut_K(Z)\, | \, \CJ (\tau ) =1\}$ where $\CJ
(\tau ) $ is the Jacobian of $\tau$ (note that
$\Aut_K(Z)=K^*\ltimes \G$ and if $K$ is {\em not} perfect then
$\im (\res) \neq \G$); $(ii)$ the bijection $\res : \Aut_K(A_1)
\ra \G$ is a monomorphism of infinite dimensional algebraic groups
which is {\em not} an isomorphism (even if $K$ is algebraically
closed); $(iii)$ an explicit formula for $\res^{-1}$ is found via
differential operators $\CD (Z)$ on $Z$ and negative powers of the
Fronenius map $F$. Proofs are based on the following (non-obvious)
equality proved in the paper:
$$ (\frac{d}{dx}+f)^p=
(\frac{d}{dx})^p+\frac{d^{p-1}f}{dx^{p-1}}+f^p, \;\; f\in K[x].$$

{\em Key Words:   the Weyl algebra, the group of automorphisms, the restriction map, the Jacobian.}

 {\em Mathematics subject classification
2000:  14J50, 16W20, 14L17,  14R10, 14R15, 14M20.}

$${\bf Contents}$$
\begin{enumerate}
\item Introduction. \item Proof of Theorem \ref{17Jul7} and the
inverse map $\th^{-1}$. \item The restriction map and its inverse.
\item The image of the restriction map $\res_n$.
\end{enumerate}
\end{abstract}


\section{Introduction}

Let $p>0$ be a prime number and $\Fp := \Z /\Z p$. Let $K$ be a
commutative $\Fp$-algebra, $A_1:=K\langle x, \der \, | \, \der x-
x\der =1\rangle$ be the first Weyl algebra over $K$. In order to
avoid awkward expressions we use $y$ instead of  $\der$ sometime,
i.e. $y=\der$. The centre $Z$ of the algebra $A_1$ is the
polynomial algebra $K[X,Y]$ in two variables $X:=x^p$ and
$Y:=\der^p$. Let $\Aut_K(A_1)$ and $\Aut_K(Z)$ be the groups of
$K$-automorphisms of the algebras $A_1$ and $Z$ respectively. They
contain the  subgroups of affine automorphisms $\Aff (A_1) \simeq
\SL_2(K)^{op}\ltimes K^2$ and $\Aff (Z)\simeq \GL_2(K)^{op}\ltimes
K^2$ respectively. If $K$ is a field of arbitrary characteristic
then the group $\Aut_K(K[X,Y])$ of automorphisms of the polynomial
algebra $K[X,Y]$ is generated by two of its subgroups, namely,
$\Aff (K[X,Y])$ and $U(K[X,Y]):=\{\phi_f : X\mapsto X$, $ Y\mapsto
Y+ f\, | \, f\in K[X]\}$. This was proved by H. W. E. Jung
\cite{Jung} in zero characteristic and by W. Van der Kulk
\cite{vdKulk} in general.

If $K$ is a field of characteristic zero J. Dixmier
\cite{DixmierBSMF68}  proved that the group $\Aut_K(A_1)$ is
generated by its subgroups $\Aff (A_1)$ and $U(A_1):=\{
\phi_{f}:x\mapsto x$, $ \der \mapsto \der + f\, | \, f\in K[x]\}$.
If $K$ is a field of characteristic $p>0$ L. Makar-Limanov
\cite{Mak-LimBSMF84} proved that the groups $\Aut_K(A_1)$ and $\G
:= \{ \tau \in \Aut_K(K[X,Y])\, | \, \CJ (\tau ) =1\}$ are
isomorphic as {\em abstract} groups where $\CJ (\tau ) $ is the
{\em Jacobian} of $\tau$. In his paper he used the restriction map
\begin{equation}\label{res1}
\res : \Aut_K(A_1) \ra \Aut_K(Z), \;\; \s \mapsto \s |_Z .
\end{equation}
In this paper, we study this map in detail. Recently, the
restriction map (for the $n$th Weyl algebra) appeared in papers of
 Y. Tsuchimoto \cite{Tsuchi-05}, A. Belov-Kanel and M. Kontsevich
 \cite{Bel-Kon05JCDP}, K.  Adjamagbo and A. van den Essen
 \cite{Adj-vdEssen06}.

\begin{theorem}\label{18Jul7}
Let $K$ be a perfect field of characteristic $p>0$. Then the
restriction map $\res$ is a group monomorphism with $\im (\res ) =
\G$.
\end{theorem}

Note that $\Aut_K(Z)= K^*\ltimes \G$ where $K^*\simeq \{ \tau_\l :
X\mapsto \l X, Y\mapsto Y\, | \, \l \in K^*\}$.

If $K$ is not perfect then Theorem \ref{18Jul7} is {\em not} true
as one can easily show  that the automorphism $\G \ni s_\mu:
X\mapsto X+\mu$, $Y\mapsto Y$, does not belong to the image of
$\res$ provided $\mu \in K\backslash F(K)$ where $F: a\mapsto a^p$
is the Frobenius. So, in the case of a perfect field we have
another proof of the result of L. Makar-Limanov (in both proofs
the results of Jung-Van der Kulk are essential).

The groups $\Aut_K(A_1)$, $\Aut_K(Z)$ and $\G$ are infinite
dimensional algebraic groups over $K$.

\begin{corollary}\label{c18Jul7}
Let $K$ be a perfect field of characteristic $p>0$. Then the
bijection  $\res : \Aut_K(A_1) \ra \G$, $\s \mapsto \s |_Z$,  is a
monomorphism of algebraic  groups over $K$ which is not an
isomorphism of algebraic groups.
\end{corollary}

The proofs of Theorem \ref{18Jul7} and Corollary \ref{c18Jul7} are
based on the following (non-obvious) formula which allows us to
find the  inverse map $\res^{-1} : \G \ra \Aut_K(A_1)$ (using
differential operators $\CD (Z)$ on $Z$), see (\ref{resm1}) and
Proposition \ref{q21Jul7}.

\begin{theorem}\label{17Jul7}
Let $K$ be a reduced commutative $\Fp$-algebra and $A_1(K)$ be the
first Weyl algebra over $K$. Then $$ (\der +f)^p=
\der^p+\frac{d^{p-1}f}{dx^{p-1}}+f^p$$ for all $f\in K[x]$. In
more detail, $(\der + f)^p= \der^p-\l_{p-1} +f^p$ where $f=
\sum_{i=0}^{p-1} \l_i x^i\in K[x]= \oplus_{i=0}^{p-1} K[x^p]x^i$,
$\l_i\in K[x^p]$.
\end{theorem}

{\it Remark}. We used the fact that $\frac{d^{p-1}f}{dx^{p-1}}=
(p-1)!\l_{p-1}$ and $(p-1)!\equiv -1\mod p$.

The group $\G$ is generated by its two subgroups $U(Z)$ and
$$ \G \cap \Aff (Z) =\{ \s_{A,a}:  \bigl(
\begin{smallmatrix} X  \\ Y
\end{smallmatrix}\bigr) \mapsto A \bigl(
\begin{smallmatrix} X  \\ Y
\end{smallmatrix}\bigr)+a \, | \, A\in \SL_2(K), a\in K^2\}\simeq \SL_2(K)^{op}\ltimes K^2.$$
Recall that the  group $\Aut_K(A_1)$ is generated by its two
subgroups $U(A_1)$ and
$$ \Aff (A_1) =\{ \s_{A,a}:  \bigl(
\begin{smallmatrix} x  \\ y
\end{smallmatrix}\bigr) \mapsto A \bigl(
\begin{smallmatrix} x  \\ y
\end{smallmatrix}\bigr)+a \, | \, A\in \SL_2(K), a\in K^2\}\simeq \SL_2(K)^{op}\ltimes K^2.$$
If $K$ is a perfect field of characteristic $p>0$ then Theorem
\ref{17Jul7} shows that
$$ \res (\Aff (A_1))= \G \cap \Aff (Z), \;\; {\rm and }\;\;  \res (U(A_1))=
U(Z).$$ In more detail, $$\res : \Aff (A_1) \ra \G \cap \Aff (Z),
\;\; \s_{\bigl(
\begin{smallmatrix} a & b\\ c & d
\end{smallmatrix}\bigr) , \bigl(
\begin{smallmatrix}  e\\  f
\end{smallmatrix}\bigr) }\mapsto \begin{cases}
\s_{\bigl(
\begin{smallmatrix} a^p & b^p\\ c^p & d^p
\end{smallmatrix}\bigr) , \bigl(
\begin{smallmatrix}  e^p\\  f^p
\end{smallmatrix}\bigr) },& \text{if $p>2$},\\
\s_{\bigl(
\begin{smallmatrix} a^2 & b^2\\ c^2 & d^2
\end{smallmatrix}\bigr) , \bigl(
\begin{smallmatrix}  e^2+ab\\  f^2+cd
\end{smallmatrix}\bigr) },& \text{if $p=2$},\\
\end{cases}
$$
see Lemma \ref{a20Jul7} and (\ref{res2}); and
$$ \res : U(A_1) \ra U(Z), \;\; \phi_f\mapsto \phi_{\th (f)},$$
where the map  $\th := F+\frac{d^{p-1}}{dx^{p-1}}:K[x]\ra K[x^p]$
is a bijection. An explicit formula for the inverse map $\th^{-1}$
is found (Proposition \ref{q21Jul7}) via differential operators
$\CD (Z)$ on $Z$ and negative powers of the Frobenius map $F$.  As
a consequence, a formula for the inverse map $\res^{-1} : \G \ra
\Aut_K(A_1)$ is given, see (\ref{resm1}).


\section{Proof of Theorem \ref{17Jul7} and the inverse map $\th^{-1}$}

In this section, a proof of Theorem \ref{17Jul7} is given and an
inversion formula for a map $\th$ is found which is a key
ingredient in the inversion formula for the restriction map.

{\bf Proof of Theorem \ref{17Jul7}}. The Weyl algebra
$A_1(K)\simeq K\t_{\Fp} A_1(\Fp )$, the Frobenius $F: a\mapsto
a^p$ and $\frac{d^{p-1}}{dx^{p-1}}$ are well behaved under ring
extensions, localizations and algebraic closure of the coefficient
 field. So, without loss of generality we may assume that $K$ is an
algebraically closed field of characteristic $p>0$: the
commutative $\Fp$-algebra $K$ is reduced, $\cap_{\gp \in \Spec
(K)}\gp =0$, and $A_1(K)/A_1(K)\gp \simeq A_1(K/ \gp )$, and so we
may assume that $K$ is a domain; then $A_1(K)\subseteq A_1({\rm
Frac} (K))\subseteq A_1(\overline{{\rm Frac} (K)})$ where  ${\rm
Frac} (K$) is the field of fractions of $K$ and $\overline{{\rm
Frac} (K)}$ is its algebraic closure.

First, let us show that the map $ L : K[x]\ra K[x^p]$, $ f\mapsto
L(f)$,  defined by the rule
$$ (\der +f)^p = \der^p + L(f) + f^p,$$
is well defined and additive, i.e. $L(f+g) = L(f) + L(g)$. The map
$$ K[x]\ra \Aut_K(A_1), f\mapsto \s_f: x\mapsto x, \der \mapsto
\der +f,$$ is a group homomorphism, i.e. $\s_{f+g} = \s_f\s_g$.
Since $\der^p\in Z(A_1) = K[x^p, \der^p]$ and $(\der +f)^p= \s
(\der )^p = \s (\der^p) \in Z(A_1)$, the map $L$ is well defined,
i.e. $L(f) \in K[x^p]$. Comparing both ends of the series of
equalities proves the additivity of the map $L$:
\begin{eqnarray*}
 \der^p+L(f+g)+f^p+g^p&=& \s_{f+g}(\der)^p = \s_{f+g} (\der^p) =
\s_f\s_g(\der^p) = \s_f(\der^p+ L(g)+g^p)\\
 &=& \der^p+ L(f)+ f^p
+L(g)+ g^p.
\end{eqnarray*}

In a view of the decomposition $K[x]=\oplus_{i=0}^{p-1}K[x^p]x^i$
and the additivity of the map $L$,  it suffices to prove the
theorem for $f=\l x^m$ where $m=0, 1, \ldots , p-1$ and $\l \in
K[x^p]$. In addition, we may assume that $\l \in K$. This follows
directly from the natural $\Fp$-algebra epimorphism
$$A_1(K[t])\ra A_1(K), \;\; t\mapsto \l ,\;\; x\mapsto x, \;\;
\der \mapsto \der, $$
 and the fact that the polynomial algebra $K[t]$ is a domain
 (hence, reduced). Therefore, it suffices to prove the theorem for
 $f=\l x^m$ where  $m=0, 1,
\ldots , p-1$ and  $\l \in K^*$.

The result is obvious for $m=0$. So, we fix the natural number $m$
such that $1\leq m \leq p-1$. Then
$$ l_m(\l ) := L(\l x^m)=\sum_{k=0}^{m-1}l_{mk}(\l )x^{kp}$$ is a sum of
{\em additive} polynomials $l_{mk}(\l )$ in $\l$ of degree $\leq
p-1$ (by the very definition of $L(\l x^m)$ and its additivity).
Recall that a polynomial $l(t) \in K[t]$ is additive if $l(\l +\mu
) = l(\l ) + \l (\mu)$ for all $\l , \mu \in K$. By Lemma 20.3.A
\cite{HumphreysLAG75}, each additive polynomial $l(t)$ is a
$p$-{\em polynomial}, i.e. a linear combination of the monomials
$t^{p^r}$, $r\geq 0$. Hence, $l_m(\l ) = a_m\l$ for some
polynomial $a_m = \sum_{k=0}^{m-1}a_{mk} x^{kp}$ where $a_{mk}\in
K$, i.e.
$$ (\der +\l x^m)^p= \der^p +\l \sum_{k=0}^{m-1}a_{mk}
x^{kp}+(\l x^m)^p.$$ Applying the $K$-automorphism $\g : x\mapsto
\mu x$, $\der \mapsto \mu^{-1}\der $, $\mu \in K^*$, of the Weyl
algebra $A_1$ to the equality above, we have

\begin{eqnarray*}
 {\rm LHS}&=& (\mu^{-1} \der + \l \mu^mx^m)^p = \mu^{-p} (\der + \l \mu^{m+1} x^m)^p \\
 &=&\mu^{-p} (\der^p +\l \mu^{m+1} \sum_{k=0}^{m-1}a_{mk}x^{kp}+(\l \mu^{m+1}x^m)^p) \\
 {\rm RHS}&=& \mu^{-p}\der^p+\l\sum_{k=0}^{m-1}a_{mk} \mu^{kp}
 x^{kp} + (\l \mu^mx^m)^p.
\end{eqnarray*}
Equating the coefficients of $x^{kp}$ gives $\l a_{mk}
\mu^{m+1-p}= \l a_{mk} \mu^{kp}$. If $a_{mk} \neq 0$ then
$\mu^{m+1-p} = \mu^{kp}$ for all $\mu \in K^*$, i.e. $m+1-p= kp$.
The maximum of $m+1-p$ is $0$ at $ m=p-1$, the minimum of $kp$ is
$0$ at $k=0$. Therefore, $a_{mk}=0$ for all $(m,k) \neq (p-1, 0)$.

For $(m,k) = (p-1, 0)$, let  $a:= a_{p-1, 0}$. Then
$$ (\der +\l x^{p-1})^p = \der^p +\l a + ( \l x^{p-1})^p.$$
In order to find the coefficient $a\in K$, consider the left
$A_1$-module
$$V:= A_1/ (A_1x^p+ A_1\der ) \simeq K[x]/ K[x^p]=
\oplus_{i=0}^{p-1}K\bx^i$$ where $\bx^i:= x^i+A_1x^p+ A_1\der$. An
easy induction on $i$ gives the equalities:
$$ (\der + \l x^{p-1})^i\bx^{p-1}= (p-1) (p-2) \cdots (p-i)
\bx^{p-1-i}, \;\; i=1,2, \ldots , p-1.$$ Now,
$$ (\der + \l x^{p-1})^p \bx^{p-1} = (\der +\l x^{p-1})(\der + \l x^{p-1})^{p-1}\bx^{p-1}=
(\der +\l x^{p-1})(p-1)!\overline{1} = (p-1)!\l \bx^{p-1}.$$ On
the other hand,
$$ (\der^p +\l a + (\l x^{p-1})^p)\bx^{p-1}= \l a \bx^{p-1}, $$
and so $a= (p-1)!\equiv -1 \mod p$. This finishes the proof of
Theorem \ref{17Jul7}. $\Box$

{\bf The map $\th$ and its inverse}. Let $K$ be a commutative
$\Fp$-algebra. The polynomial algebra $K[x]= \oplus_{i\geq 0}
Kx^i$ is a positively graded algebra and a positively filtered
algebra $K[x] = \cup_{i\geq 0} K[x]_{\leq i}$ where $K[x]_{\leq
i}:= \oplus_{j=0}^iKx^j= \{ f\in K[x]\, | \, \deg (f) \leq i\}$.
Similarly, the polynomial algebra $K[x^p]$ in the variable $x^p$
is a positively  graded algebra $K[x^p]= \oplus_{i\geq 0} Kx^{pi}$
 and a positively filtered algebra
$K[x^p] = \cup_{i\geq 0} K[x^p]_{\leq i}$ where $K[x^p]_{\leq i}:=
\oplus_{j=0}^iKx^{pj}= \{ f\in K[x^p]\, | \, \deg_{x^p} (f) \leq
i\}$. The associated graded algebras ${\rm gr} \, K[x]$ and ${\rm
gr} \, K[x^p]$ are canonically isomorphic to $K[x]$ and $K[x^p]$
respectively.  For a polynomial $f=\sum_{i=0}^d\l_ix^i\in K[x]$
(resp. $g=\sum_{i=0}^d\mu_ix^{pi}\in K[x^p]$) of degree $d$,
$\l_dx^d$ (resp. $\mu_dx^{pd}$) is called the leading term of $f$
(resp. $g$) denoted $l(f)$ (resp. $l(g)$). Consider the
$\Fp$-linear map (see Theorem \ref{17Jul7}) 
\begin{equation}\label{thFp}
\th : F+\frac{d^{p-1}}{dx^{p-1}}: K[x]\ra K[x^p], f\mapsto
f^p+\frac{d^{p-1}f}{dx^{p-1}},
\end{equation}
where $F: f\mapsto f^p$ is the Frobenius ($\Fp$-algebra
monomorphism). In more detail, $$\th : K[x]= \oplus_{i=0}^{p-1}
K[x^p]x^i\ra K[x^p]= \oplus_{i=0}^{p-1} K[x^{p^2}]x^{pi}, \;\;
\sum_{i=0}^{p-1} a_ix^i\mapsto \sum_{i=0}^{p-1}
a_i^px^{pi}-a_{p-1}, $$ where $a_i\in K[x^p]$. This means that the
map  $\th$ respects the filtrations of the algebras $K[x]$ and
$K[x^p]$, $\th (K[x]_{\leq j})\subseteq K[x^p]_{\leq j}$ for all
$j\geq 0$, and so the associated graded map ${\rm gr} (\th
):K[x]\ra K[x^p]$ coincides with the Frobenius $F$,
\begin{equation}\label{grt=F}
{\rm gr} (\th ) = F.
\end{equation}

\begin{lemma}\label{p21Jul7}
Let $K$ be a perfect field of characteristic $p>0$. Then
\begin{enumerate}
\item ${\rm gr} (\th ) = F: K[x]\ra K[x^p]$ is an isomorphism of
$\Fp$-algebras. \item  $\th : K[x]\ra K[x^p]$ is the isomorphism
of vector spaces over $\Fp$ such that $\th ( K[x]_{\leq i})=
K[x^p]_{\leq i}$, $i\geq 0$. \item For each $f\in K[x]$, $l(\th
(f))= l(f)^p$.
\end{enumerate}
\end{lemma}

{\it Proof}.  Statement 1 is obvious since $K$ is a perfect field
of characteristic $p>0$ ($F(K)=K$). Statements 2 and 3 follow from
statement 1. $\Box $

{\it Remark}. The problem of finding the inverse map $\res^{-1}$
of the group isomorphism $\res : \Aut_K(A_1)\ra \G$, $\s \mapsto
\s |_Z$, is essentially equivalent to the problem of finding
$\th^{-1}$, see (\ref{resm1}).

The inversion formula for $\th^{-1}$ (Proposition \ref{q21Jul7})
is given via certain differential operators. We recall some facts
on differential operators that are needed in the proof of
Proposition \ref{q21Jul7}.

Let $K$ be a field of characteristic $p>0$ and $\CD (K[x]) =
\oplus_{i\geq 0}K[x]\der^{[i]}$ be the ring of differential
operators on the polynomial algebra $K[x]$ where $\der^{[i]} :=
\frac{\der^i}{i!}$. The algebra $K[x]$ is the left $\CD
(K[x])$-module (in the usual sense). The subalgebra $K[x^p]=\oplus
_{i=0}^{p-1} K[x^{p^2}]x^{pi}$ of $K[x]$ is
$x^p\der^{[p]}$-invariant and, for each $i=0, 1, \ldots , p-1$, $
K[x^{p^2}]x^{pi}$ is the eigenspace of the element $x^p\der^{[p]}$
that corresponds to the eigenvalue $i$. Let $J(i) := \{ 0,1,
\ldots , p-1\} \backslash \{ i \}$. Then 
\begin{equation}\label{projpi}
\pi_i:= \der^{[pi]}\frac{\prod_{j\in
J(i)}(x^p\der^{[p]}-j)}{\prod_{j\in J(i)}(i-j)}: K[x^p]\ra
K[x^{p^2}], \;\; \sum_{i=0}^{p-1}a_i x^{pi}\mapsto a_i,
\end{equation}
where all $a_i\in K[x^{p^2}]$ (since the map $\frac{\prod_{j\in
J(i)}(x^p\der^{[p]}-j)}{\prod_{j\in J(i)}(i-j)}: K[x^p]\ra K[x^p]$
is the projection onto the summand $K[x^{p^2}]x^{pi}$ in the
decomposition $K[x]=\oplus_{i=0}^{p-1} K[x^{p^2}]x^{pi}$ and
$\der^{[pi]}(a_ix^{pi})= a_i$).

Let $K$ be a perfect field of characteristic $p>0$. Consider the
$\Fp$-linear map 
\begin{equation}\label{dp1F1}
\der^{[(p-1)p]}F^{-1}:K[x^{p^2}]\ra K[x^{p^2}], \;\; \sum_{i\geq
0}a_ix^{p^2i}\mapsto \sum_{i\geq 0}
a^\frac{1}{p}_{p-1+pi}x^{p^2i},
\end{equation}
where $a_i\in K$. By induction on a natural number $n$, we have
\begin{equation}\label{dp1F2}
(\der^{[(p-1)p]}F^{-1})^n (\sum_{i\geq 0}a_ix^{p^2i})=\sum_{i\geq
0} a^{p^{-n}}_{(p-1)(1+p+\cdots + p^{n-1})+p^ni}x^{p^2i}, \;\;
n\geq 1.
\end{equation}
This shows that the map $\der^{[(p-1)p]}F^{-1}$ is  a {\em locally
nilpotent map}. This means that  $K[x^{p^2}]= \cup_{n\geq 1} \ker
(\der^{[(p-1)p]}F^{-1})^n$, i.e. for each element $a\in
K[x^{p^2}]$, $(\der^{[(p-1)p]}F^{-1})^n(a)=0$ for all $n\gg 0$.
Hence, the map $1-\der^{[(p-1)p]}F^{-1}$ is invertible and its
inverse is given by the rule

\begin{equation}\label{dp1F3}
(1-\der^{[(p-1)p]}F^{-1})^{-1}=\sum_{j\geq 0}
(\der^{[(p-1)p]}F^{-1})^j.
\end{equation}

The following proposition gives an explicit formula for
$\th^{-1}$.
\begin{proposition}\label{q21Jul7}
Let $K$ be a perfect field of characteristic $p>0$. Then the
inverse map $\th^{-1} : K[x^p]=\oplus_{i=0}^{p-1}
K[x^{p^2}]x^{pi}\ra K[x]=\oplus _{i=0}^{p-1} K[x^p]x^i$,
$\sum_{i=0}^{p-1}\mu_ix^{pi}\mapsto \sum_{i=0}^{p-1} \l_ix^i$,
$\mu_i\in K[x^{p^2}]$, $\l_i\in K[x^p]$, is given by the rule
\begin{enumerate}
\item for $i=0,1, \ldots , p-2$, $\l_i= \mu_i^\frac{1}{p}+F^{-1}
\pi_iF^{-1} \sum_{j\geq 0} (\der^{[(p-1)p]}F^{-1})^j(\mu_{p-1})$,
\item $\l_{p-1}=(\sum_{i=0}^{p-2}x^{pi} \pi_i F^{-1} \sum_{j\geq
0} (\der^{[(p-1)p]}F^{-1})^j +x^{p(p-1)}\sum_{j\geq
1}(\der^{[(p-1)p]}F^{-1})^j)(\mu_{p-1})$ where $\pi_i$ is defined
in (\ref{projpi}).
\end{enumerate}
\end{proposition}

{\it Proof}. Let $g=\sum_{i=0}^{p-1}\mu_ix^{pi}\in K[x^p]$,
$\mu_i\in K[x^{p^2}]$; $f=\sum_{i=0}^{p-1} \l_ix^i\in K[x]$,
$\l_i\in K[x^p]$; and $\l_{p-1} = \sum_{i=0}^{p-1} a_i x^{pi}$,
$a_i\in K[x^{p^2}]$. Then $\th^{-1} (g) = f$ iff $g= \th (f)$ iff
$F^{-1} (g) = F^{-1} \th (f)$ iff
$$ \sum_{i=0}^{p-1} F^{-1} (\mu_i) x^i = F^{-1} (F(f)-\l_{p-1}) =
f-F^{-1} (\l_{p-1}) = \sum_{i=0}^{p-1} (\l_i-F^{-1} (a_i))x^i$$
iff

\begin{equation}\label{ith1}
\l_i= F^{-1} (\mu_i+a_i), \;\; i=0, 1, \ldots , p-1.
\end{equation}
For $i=p-1$, the equality (\ref{ith1}) can be rewritten as follows
\begin{equation}\label{ith2}
\sum_{i=0}^{p-2}a_ix^{pi}+a_{p-1}x^{p(p-1)}= F^{-1}
(\mu_{p-1}+a_{p-1}).
\end{equation}
For each $i=0, 1, \ldots , p-2$, applying the map  $\pi_i$ (see
(\ref{projpi})) to (\ref{ith2}) gives the equality $a_i=
\pi_iF^{-1} (\mu_{p-1}+a_{p-1})$, and so the equalities
(\ref{ith1}) can be rewritten as follows

\begin{equation}\label{ith3}
\l_i= F^{-1}(\mu_i+\pi_iF^{-1}(\mu_{p-1}+a_{p-1})), \;\; i=0,1,
\ldots , p-2.
\end{equation}
Applying $\der^{[(p-1)p]}$  to (\ref{ith2}) yields $a_{p-1} =
\der^{[(p-1)p]}F^{-1}(\mu_{p-1}+a_{p-1})$, and so $(1-\D )a_{p-1}=
\D (\mu_{p-1})$ where $\D : = \der^{[(p-1)p]}F^{-1}$. By
(\ref{dp1F3}), $a_{p-1}=\sum_{j\geq 1}\D^j(\mu_{p-1})$. Putting
this expression in (\ref{ith3}) yields,  $$\l_i = F^{-1} (\mu_i)
+F^{-1} \pi_i  F^{-1} \sum_{j\geq 0} \D^j (\mu_{p-1}), \;\; i=0,1,
\ldots , p-2.$$ This proves statement 1.  Finally,
\begin{eqnarray*}
 \l_{p-1} &= & \sum_{i=0}^{p-1}a_ix^{pi}= \sum_{i=0}^{p-2} a_ix^{pi} + a_{p-1}x^{p(p-1)}
 \\ &=& \sum_{i=0}^{p-2} x^{pi} \pi_i F^{-1} (\mu_{p-1}+a_{p-1}) + x^{p(p-1)}\sum_{j\geq 1} \D^j(\mu_{p-1}) \\
 &=&\sum_{i=0}^{p-2} x^{pi} \pi_i F^{-1} \sum_{j\geq 0} \D^j
 (\mu_{p-1})+x^{p(p-1)}\sum_{j\geq 1} \D^j(\mu_{p-1}) \\
 &=&(\sum_{i=0}^{p-2}x^{pi} \pi_i F^{-1} \sum_{j\geq
0} (\der^{[(p-1)p]}F^{-1})^j +x^{p(p-1)}\sum_{j\geq
1}(\der^{[(p-1)p]}F^{-1})^j)(\mu_{p-1}). \;\;\; \Box
\end{eqnarray*}


\section{The restriction map and its inverse}

In this section, Theorem \ref{18Jul7}, \ref{19Jul7} and Corollary
\ref{c18Jul7} are proved. An inversion formula for the restriction
map $\res : \Aut_K(A_1) \ra \G$ is found, see (\ref{resm1}).

{\bf The group of affine automorphisms}. Let $K$ be a perfect
field of characteristic $p>0$. Each element $a$ of the Weyl
algebra $A_1= \oplus_{i,j\in \N}Kx^iy^i$ is a unique sum $a= \sum
\l_{ij} x^iy^j$ where all but finitely many scalars $\l_{ij}\in K$
are equal to zero. The number $\deg (a):=\max \{ i+j\, | \,
\l_{ij}\neq 0\}$ is called the degree of $a$, $\deg (0):=-\infty$.
Note that $\deg (ab) = \deg (a) +\deg (b)$, $\deg (a+b)\leq \max
\{ \deg (a) , \deg (b)\}$, and $\deg (\l a) = \deg (a)$ for all
$\l \in K^*$. For each $\s \in \Aut_K(A_1)$,
$$\deg (\s ) := \max
\{ \deg (\s (x)), \deg (\s (y))\}$$
 is called the {\em degree} of
$\s$. The set (which is obviously a subgroup of $\Aut_K(A_1)$)
$\Aff (A_1) =\{ \s \in \Aut_K(A_1)\, | \, \deg (\s ) = 1\}$ is
called the group of affine automorphisms of the Weyl algebra
$A_1$. Clearly,
$$\Aff (A_1) = \{ \s_{A, a} : \bigl(
\begin{smallmatrix} x  \\ y
\end{smallmatrix}\bigr) \mapsto A \bigl(
\begin{smallmatrix} x  \\ y
\end{smallmatrix}\bigr)+a \, | \, A\in \SL_2(K), a\in K^2\}, \;\; \s_{A,a}\s_{B,b}= \s_{BA, Ba+b}.$$
For each group $G$, let $G^{op}$ be its {\em opposite} group
($G^{op}=G$ as sets but the product $ab$ in $G^{op}$ is equal to
 $ba$ in $G$). The map $G\ra G^{op}$, $g\mapsto g^{-1}$, is a
 group
automorphism. The group $\Aff (A_1)$ is the semi-direct product
$\SL_2(K)^{op}\ltimes K^2$ of its subgroups $\SL_2(K)^{op} = \{
\s_{A, 0} \, | \, A\in \SL_2(K)\}$ and $K^2\simeq  \{ \s_{1, a}\,
| \, a\in K^2\}$ where $K^2$ is the normal subgroup of $\Aff
(A_1)$  since $\s_{A,0}\s_{1,a}\s_{A,0}^{-1} = \s_{1, A^{-1} a}$.
It is obvious that the group $\Aff (A_1)$ is generated by the
automorphisms:
$$ s: x\mapsto y, \; y\mapsto -x; \;\; t_\mu : x\mapsto \mu x,\;
y\mapsto \mu^{-1} y; \;\; \phi_{\l x^i} : x\mapsto x, \; y\mapsto
y +\l x^i,
$$
where $\l \in K$, $\mu \in K^*$ and $i=0, 1$.

Recall that the centre $Z$ of the Weyl algebra $A_1$ is the
polynomial algebra $K[X,Y]$ in $X:= x^p$ and $Y:= y^p$ variables.
Let $\deg (z )$ be the total degree of a polynomial $z\in Z$. For
each automorphisms $\s \in \Aut_K(Z)$,
$$\deg (\s ) := \max \{ \deg
(\s (X)), \deg (\s (Y))\}$$ is called the {\em degree} of $\s$.
$$\Aff (Z) := \{ \s \in \Aut_K(Z)\, | \, \deg (\s ) =1\} =
\{ \s_{A, a} : \bigl(
\begin{smallmatrix} x  \\ y
\end{smallmatrix}\bigr) \mapsto A \bigl(
\begin{smallmatrix} x  \\ y
\end{smallmatrix}\bigr)+a \, | \, A\in \GL_2(K), a\in K^2\}$$
is the group of affine automorphisms of $Z$, $ \s_{A,a}\s_{B,b}=
\s_{BA, Ba+b}$, The group $\Aff (A_1)$ is the semi-direct product
$\GL_2(K)^{op}\ltimes K^2$ of its subgroups $\GL_2(K)^{op} = \{
\s_{A, 0} \, | \, A\in \GL_2(K)\}$ and $K^2\simeq  \{ \s_{1, a}\,
| \, a\in K^2\}$ where $K^2$ is the normal subgroup of $\Aff (Z)$
since $\s_{A,0}\s_{1,a}\s_{A,0}^{-1} = \s_{1, A^{-1} a}$.

A group $G$ is called an {\em exact} product of its subgroups
$G_1$ and $G_2$ denoted $ G= G_1\times_{ex} G_2$ if each element
$g\in G$ is a unique product $g=g_1g_2$ for some elements $g_1\in
G_1$ and $g_2\in G_2$. Then $ \GL_2(K)^{op} = K^*\times_{ex}
\SL_2(K)^{op}$ where $K^*\simeq \{ \g_\mu : X\mapsto \mu X, \;
Y\mapsto Y\, | \, \mu \in K^*\}$, $\g_\mu \g_\nu = \g_{\mu \nu}$.
Clearly, $\Aff (Z)= (K^*\times_{ex} \SL_2(K)^{op})\ltimes K^2$,
and so the group $\Aff (Z)$ is generated by the following
automorphisms (where $\l \in K$, $\mu \in K^*$ and $i=0,1$):
$$ s: X\mapsto Y, \; Y\mapsto -X; \;\; t_\mu : X\mapsto \mu X,\;
Y\mapsto \mu^{-1} Y; \phi_{\l X^i} : X\mapsto X, \; Y\mapsto Y +\l
X^i; \;\; {\rm and } \;\; \g_\mu .
$$
The automorphisms $t_\mu$ and $\g_\nu$ commute.
\begin{lemma}\label{a20Jul7}
Let $K$ be a perfect field of characteristic $p>0$. Then the
restriction map $\res_{aff}:\Aff (A_1) \ra \Aff (Z)$, $\s \mapsto
\s |_Z$, is a group monomorphism with $\im (\res_{aff})=
\SL_2(K)^{op} \ltimes K^2$.
\end{lemma}

{\it Proof}. Since $\res_{aff}(s) = s$; $\res_{aff}(t_\mu )=
t_{\mu^p}$; for $i=0,1$,   $\res_{aff} (\phi_{\l x^i})= \phi_{\l^p
X^i}$ if $p>2$ and  $\res_{aff} (\phi_{\l x^i})= \phi_{\l^2
X^i+\d_{i,1}\l}$ if $p=2$ where $\d_{i,1}$ is the Kronecker delta
(Theorem \ref{17Jul7}), i.e. 
\begin{equation}\label{res2}
\res_{aff} (\s_{\bigl(
\begin{smallmatrix} a & b\\ c & d
\end{smallmatrix}\bigr) , \bigl(
\begin{smallmatrix}  e\\  f
\end{smallmatrix}\bigr) })=\begin{cases}
\s_{\bigl(
\begin{smallmatrix} a^p & b^p\\ c^p & d^p
\end{smallmatrix}\bigr) , \bigl(
\begin{smallmatrix}  e^p\\  f^p
\end{smallmatrix}\bigr) },& \text{if $p>2$},\\
\s_{\bigl(
\begin{smallmatrix} a^2 & b^2\\ c^2 & d^2
\end{smallmatrix}\bigr) , \bigl(
\begin{smallmatrix}  e^2+ab\\  f^2+cd
\end{smallmatrix}\bigr) },& \text{if $p=2$},\\
\end{cases}
\end{equation}
the result is obvious.  $\Box $

\begin{lemma}\label{b20Jul7}
The automorphisms of the algebra $Z$:  $s$, $t_\mu$,  $\phi_{\l
X^i}$ and $\g_\mu$ satisfy the following relations:
\begin{enumerate}
\item  $st_\mu = t_{\mu^{-1}}s$ and $s\g_\mu = \g_\mu
t_{\mu^{-1}}s$. \item $\phi_{\l X^i}t_\mu = t_\mu \phi_{\l
\mu^{-i-1} X^i}$ and $\phi_{\l X^i} \g_\mu = \g_\mu \phi_{\l
\mu^{-i} X^i}$. \item $s^2= t_{-1}$, $s^{-1} = t_{-1}s: X\mapsto
-Y$, $Y\mapsto X$.
\end{enumerate}
\end{lemma}

{\it Proof}.  Straightforward. $\Box $

The map
$$ K[X]\ra \Aut(Z), \;\; f\mapsto \phi_f: X\mapsto X, \;\;
Y\mapsto Y+f,$$ is a group monomorphism ($\phi_{f+g} =
\phi_f\phi_g$). For $\s \in \Aut (Z)$, $\CJ (\s ) := \det \bigl(
\begin{smallmatrix} \frac{\der \s (X)}{\der X} & \frac{\der \s (X)}{\der Y }\\ \frac{\der \s (Y)}{\der X}
&\frac{\der \s (Y)}{\der Y}
\end{smallmatrix}\bigr)$ is the Jacobian of $\s$. It follows from
the equality (which is a direct consequence of the chain rule)
$\CJ ( \s \tau ) = \CJ ( \s ) \s (\CJ (\tau ))$ that $\CJ (\s )
\in K^*$ (since $1= \CJ (\s \s^{-1} ) = \CJ (\s ) \s (\CJ
(\s^{-1}))$ in $K[X,Y]$), and so the kernel $\G := \{ \s \in
\Aut_K(Z)\, | \, \CJ (\s ) = 1\}$ of the group epimorphism $\CJ :
\Aut (Z)\ra K^*$, $\s \mapsto \CJ (\s )$, is a normal subgroup of
$\Aut_K(Z)$. Hence,  
\begin{equation}\label{}
\Aut_K(Z)= K^*\ltimes \G
\end{equation}
 is the semi-direct product of its subgroups $\G$ and $K^*\simeq
 \{ \g_\mu \, | \, \mu \in K^*\}$.
\begin{corollary}\label{c20Jul7}
Let $K$ be a field of characteristic $p>0$. Then
\begin{enumerate}
\item Each automorphism $\s \in \Aut_K(Z)$ is a product $\s =
\g_\mu t_\nu \phi_{f_1}s\phi_{f_2}s\ldots
\phi_{f_{n-1}}s\phi_{f_n}$ for some $\mu , \nu \in K^*$ and
$f_i\in K[x]$. \item Each automorphism $\s \in \G$ is a product
$\s =t_\nu \phi_{f_1}s\phi_{f_2}s\ldots \phi_{f_{n-1}}s\phi_{f_n}$
for some $\nu \in K^*$ and $f_i\in K[x]$.
\end{enumerate}
\end{corollary}

{\it Proof}. 1. Statement 1 follows at one from Lemma
\ref{b20Jul7} and the fact that the group $\Aut_K(Z)$ is generated
by $\Aff (Z)$ and $\phi_{\l X^i}$, $\l \in K$, $i\in \N$.

2. Statement 2 follows from statement 1: $\s =\g_\mu t_\nu
\phi_{f_1}s\phi_{f_2}s\ldots \phi_{f_{n-1}}s\phi_{f_n}\in \G$ iff
$$ 1= \CJ (\s ) = \CJ (\g_\mu t_\nu
\phi_{f_1}s\phi_{f_2}s\ldots \phi_{f_{n-1}}s\phi_{f_n})=\CJ
(\g_\mu ) \g_\mu (1) = \mu $$ iff $\s =t_\nu
\phi_{f_1}s\phi_{f_2}s\ldots \phi_{f_{n-1}}s\phi_{f_n}$.  $\Box $

{\bf Proof of Theorem \ref{18Jul7}}. {\em Step 1:  $\res$ is a
monomorphism}. It is obvious that 
\begin{equation}\label{res3}
\deg \, \res (\s ) = \deg \, \s, \;\; \s \in \Aut_K(A_1).
\end{equation}
The map $\res$ is a group homomorphism, so we have to show that
$\res ( \s ) = {\rm id}_Z$ implies $\s = {\rm id}_{A_1}$ where
${\rm id}_Z$ and ${\rm id}_{A_1}$ are the identity maps on $Z$ and
$A_1$ respectively. By (\ref{res3}), $\res ( \s ) = {\rm id}_Z$
implies $\deg (\s ) = 1$. Then, by (\ref{res2}), $\s = {\rm
id}_{A_1}$.

 {\em Step 2:  $\G \subseteq \im (\res )$}. By Corollary
 \ref{c20Jul7}.(2), each automorphism $\s \in \G$ is a product
 $\s =t_\nu \phi_{f_1}s\ldots \phi_{f_{n-1}}s\phi_{f_n}$. Since
 $\res (t_{\nu^\frac{1}{p}})= t_\nu$, $\res (\phi_{\th^{-1}
 (f_i)})= \phi_{f_i}$ and $\res (s) = s$, we  have
 $\s =\res (t_{\nu^\frac{1}{p}} \phi_{\th^{-1}(f_1)}s\ldots
 \phi_{\th^{-1}(f_{n-1})}s\phi_{\th^{-1}(f_n)})$, and so
 $\G \subseteq \im (\res )$.

{\em Step 3:  $\G = \im (\res )$}. Let $\s \in \im (\res )$. By
Corollary \ref{c20Jul7}.(1), $$\res ( \s ) =\g_\mu t_\nu
\phi_{f_1}s\ldots \phi_{f_{n-1}}s\phi_{f_n}=\g_\mu \res (\tau )$$
for some $\tau \in \Aut_K(A_1)$ such that $\res (\tau ) \in \G$,
by Step 2. Then $\res (\s \tau^{-1}) = \g_\mu$. By (\ref{res3}),
$\deg (\s \tau^{-1}) = \deg \, \res ( \s \tau^{-1})= \deg \,
\g_\mu = 1$, and so $\s\tau^{-1} \in \Aff (A_1)$. By Lemma
\ref{a20Jul7}, $\g_\mu =1$, and so $\s = \tau$, hence $\res (\s )
= \res (\tau )\in \G$. This means that $\G = \im (\res )$. $\Box$

If $K$ is a {\em perfect} field of characteristic $p>0$ we obtain
the result of L. Makar-Limanov.

\begin{theorem}\label{19Jul7}
Let $K$ be a perfect field of characteristic $p>0$. Then the group
$\Aut_K(A_1)$ is generated by $\Aff (A_1)\simeq
\SL_2(K)^{op}\ltimes K^2$ and the automorphisms $\phi_{\l x^i}$,
$\l \in K^*$, $i=2,3, \ldots $.
\end{theorem}

{\it Proof}.  By Theorem \ref{18Jul7}, the map $\res :\Aut_K(A_1)
\ra \G$ is the isomorphism of groups. By Corollary
\ref{c20Jul7}.(2), each element $\g \in \G$ is a product $$\g =
t_\nu \phi_{f_1} s\ldots \phi_{f_{n-1}}s\phi_{f_n}=\res
(t_{\nu^\frac{1}{p}} \phi_{\th^{-1}(f_1)}s\ldots
 \phi_{\th^{-1}(f_{n-1})}s\phi_{\th^{-1}(f_n)}).$$
Now, it is obvious that the group $\Aut _K(A_1)$ is generated by
$\Aff (A_1)$ and the automorphisms $\phi_{\l x^i}$, $\l \in K^*$,
$i=2,3, \ldots $.  $\Box $

{\bf The inverse map $\res^{-1} : \G \ra \Aut_K(A_1)$}. By
Corollary \ref{c20Jul7}.(2), each element $\g \in \G$ is a product
 $\g = t_\nu \phi_{f_1} s\ldots \phi_{f_{n-1}}s\phi_{f_n}$. By
 Proposition \ref{q21Jul7}, the inverse map for $\res$ is given by
 the rule

\begin{equation}\label{resm1}
\res^{-1}:\G \ra \Aut_K(A_1), \; \g = t_\nu \phi_{f_1} s\ldots
\phi_{f_{n-1}}s\phi_{f_n}\mapsto t_{\nu^\frac{1}{p}}
\phi_{\th^{-1}(f_1)}s\ldots
 \phi_{\th^{-1}(f_{n-1})}s\phi_{\th^{-1}(f_n)}.
\end{equation}

{\bf Proof of Corollary \ref{c18Jul7}}. The group $\Aut_K(A_1)$
(resp. $\Aut_K(Z)$) are infinite-dimensional algebraic groups over
$K$ (and over $\Fp$) where the coefficients of the polynomials $\s
(x)$ and $\s (y)$ where $\s \in \Aut_K(A_1)$ (resp. of $\tau (X)$
and $\tau (Y)$ where $\tau \in \Aut_K(Z)$) are coordinate
functions, i.e. generators for the algebra of regular functions on
$\Aut_K(A_1)$ (resp. $\Aut_K(Z)$). The group $\G$ is a closed
subgroup of $\Aut_K(Z)$. By the very definition, the map $\res :
\Aut_K(A_1) \ra \G$ is a polynomial map (i.e. a morphism of
algebraic varieties). By (\ref{resm1}) and Proposition
\ref{q21Jul7}, $\res^{-1}$ is not a polynomial map over (and over
$\Fp$ as well). $\Box$


\section{The image of the
restriction map $\res_n$}

Let $K$ be a field of characteristic $p>0$ and $A_n= K\langle x_1,
\ldots , x_{2n}\rangle$ be the $n$th Weyl algebra over $K$: for
$i,j=1, \ldots , n$,
$$ [x_i, x_j]=0, \;\;  [x_{n+i}, x_{n+j}]=0, \;\; [x_{n+i},
x_i]=\d_{ij},$$ where $\d_{ij}$ is the Kronecker delta. The centre
$Z_n$ of the algebra $A_n$ is the polynomial algebra $K[X_1,
\ldots , X_{2n}]$ in $2n$ variables where $X_i:=x_i^p$. The groups
of $K$-automorphisms $\Aut_K(A_n)$ and $\Aut_K(Z_n)$  contain the
affine subgroups $\Aff (A_n) = {\rm Sp}_{2n}(K)^{op}\ltimes K^n$
and $\Aff (Z_n) = \GL_n(K)^{op}\ltimes K^n$ respectively. Clearly,
 $\Aff (A_n) = \{ \s \in \Aut_K(A_n) \, | \, \deg (\s ) = 1\}$ and
 $\Aff (Z_n)= \{ \tau \in \Aut_K(Z_n)\, | \, \deg (\tau ) = 1\}$
 where $\deg (\s )$ (resp. $\deg (\tau )$) is the (total)  degree of $\s $
 (resp. $\tau $) defined in the obvious way. The kernel $\G_n$ of
 the group epimorphism $\CJ : \Aut_K(Z_n) \ra K^*$, $\tau \mapsto
 \CJ (\tau ) := \det ( \frac{\der \tau (X_i)}{\der X_j})$ is the
 normal subgroup $\G_n := \{ \tau \in \Aut_K(Z_n) \, | \, \CJ
 (\tau ) =1\}$, and $\Aut_K(Z_n) = K^*\ltimes \G_n$ is the
 semi-direct product of $K^*\simeq \{ \g_\mu  \, | \, \g_\mu (X_1) =
 \mu X_1, \g_\mu (X_j) = X_j, j=2, \ldots , 2n;  \mu \in K^*\}$ and $\G_n$.

By considering leading terms of the polynomials $\s (X_i)$, it
follows as in the case $n=1$ that the restriction map
$$ \res_n : \Aut_K(A_n) \ra \Aut_K(Z_n), \;\; \s \mapsto \s
|_{Z_n}, $$ is a group monomorphism. If $K$ is a perfect field
then
$$ \res_n (\Aff (A_n)) = {\rm Sp}_{2n}(K)^{op}\ltimes
K^{2n}\subset \Aff (Z_n) = \GL_{2n}(K)^{op}\ltimes K^{2n}.$$ This
follows from the fact that for any element of $\Aff (A_n)$,
$\s_{A, a}:x\mapsto Ax+a$ where $A=(a_{ij})\in {\rm Sp}_{2n}(K)$
and $a= (a_i)\in K^{2n}$, 
\begin{equation}\label{resnAa}
\res_n(\s_{A,a})=\begin{cases}
\s_{(a_{ij}^p), (a_i^p)}& \text{if $p>2$},\\
\s_{(a_{ij}^2), (a_i^2+\sum_{j=1}^na_{ij} a_{i,n+j})}& \text{if $p=2$},\\
\end{cases}
\end{equation}
which can be proved in the same fashion as (\ref{res2}). Since
${\rm Sp}_{2n}(K)\subseteq \SL_{2n}(K)$ (any symplectic matrix
$S\in {\rm Sp}_{2n}(K)$ has the from $S= TJT^{-1}$ for some matrix
$T\in \GL_{2n}(K)$ where $J= {\rm diag} (\bigl(
\begin{smallmatrix}0  & 1\\ -1 & 0
\end{smallmatrix}\bigr), \ldots , \bigl(
\begin{smallmatrix}0  & 1\\ -1 & 0
\end{smallmatrix}\bigr))$, $n$ times, hence $\det (S)= 1$),
$$ \res_n (\Aff (A_n))\subseteq \SL_{2n}(K)^{op}\ltimes
K^{2n}\subset \G_n.$$

{\it Question 1. For an algebraically closed field $K$ of
characteristic $p>0$, is $\im (\res_n)\subseteq \G_n$?}

{\it Question 2. For an algebraically closed field $K$ of
characteristic $p>0$, is the injection}
$$\Aff (Z_n) / \res_n (\Aff (A_n))\simeq \GL_{2n}(K)^{op}/{\rm
Sp}_{2n}(K)^{op}\ra \Aut_K(Z_n)/ \im (\res_n)$$ {\it a bijection?}

\begin{corollary}\label{22Jul7}
Let $K$ be a reduced commutative $\Fp$-algebra, $A_n(K)$ be the
Weyl algebra, and $\der_i:= x_{n+i}$. Then
$$ (\der_i +f)^p= \der_i^p+\frac{\der^{p-1}f}{\der x_i^{p-1}}+f^p$$
 for all polynomials $f\in K[x_1, \ldots , x_n]$.
\end{corollary}

{\it Proof}. Without loss of generality we may assume that $i=1$.
Since $K[x_2, \ldots , x_n]$ is a reduced commutative
$\Fp$-algebra and $\der_1+f\in A_1(K[x_2, \ldots , x_n])$, the
result follows from Theorem \ref{17Jul7}.  $\Box $

Department of Pure Mathematics

University of Sheffield

Hicks Building

Sheffield S3 7RH

UK

email: v.bavula@sheffield.ac.uk

\end{document}